\documentclass[article]{elsarticle}

\usepackage[utf8]{inputenc}
\usepackage{amssymb}
\usepackage{amsfonts}
\usepackage{graphicx}
\usepackage{multirow}
\usepackage{bigstrut}
\usepackage{epstopdf}
\usepackage[english]{babel}
\usepackage{combelow}
\usepackage{amsmath}
\usepackage{amsthm}
\usepackage{algorithm}
\usepackage{algorithmicx}
\usepackage[noend]{algpseudocode}
\usepackage{xcolor}
\usepackage{bm}

\usepackage{tikzpagenodes}
\usepackage{tikz}
\usetikzlibrary{shapes}
\usetikzlibrary{arrows}
\usetikzlibrary{positioning}
\usetikzlibrary{matrix}
\usetikzlibrary{patterns}
\usetikzlibrary{decorations.pathreplacing,decorations.pathmorphing,decorations.markings}
\usepackage{tikz-qtree}
\usetikzlibrary{trees,calc,arrows.meta,positioning,bending}
\usepackage{mathtools}
\usepackage{todonotes}
\usepackage{sidecap}
\usepackage{fullpage}
\usepackage{hyperref}
\usepackage{upgreek}

\newtheorem{remark}{Remark}

\graphicspath{{figures/}}
\DeclareGraphicsExtensions{.pdf,.eps,.png,.jpg,.jpeg}

\newtheorem{lemma}{Lemma}

\hypersetup{
    pdftitle={Distributed computing for physics-based data-driven reduced modeling at scale: Application to a rotating detonation rocket engine},
    pdfauthor={Ionut-Gabriel Farcas, Rayomand Gundevia, Ramakanth Munipalli, Karen Willcox},
    pdfkeywords={distributed computing, distributed-memory computing, distributed-memory data-driven learning, distributed computing for data-driven learning, physics-based data-driven learning, physics-based data-driven reduced modeling, data-driven reduced modeling, scientific machine learning, large-scale simulation, rotating detonation rocket engines, combustion}
}
\usepackage{fullpage}
\usepackage{hyperref}

\journal{Computer Physics Communications}

\begin{document}

\begin{frontmatter}

\title{Distributed computing for physics-based data-driven reduced modeling at scale: Application to a rotating detonation rocket engine}

\author[label1,label2]{Ionu\cb{t}-Gabriel Farca\cb{s}}
\address[label1]{Department of Mathematics, Virginia Tech, Blacksburg VA, USA}
\ead{farcasi@vt.edu}
\address[label2]{Oden Institute for Computational Engineering and Sciences, The University of Texas at Austin, Austin TX, USA}
\ead{ionut.farcas@austin.utexas.edu}

\author[label3]{Rayomand P. Gundevia}
\address[label3]{Amentum, Edwards Air Force Base, Edwards CA, USA}
\ead{rayomand.gundevia.ctr@afrl.af.mil}

\author[label4]{Ramakanth Munipalli}
\address[label4]{Air Force Research Laboratory, Edwards Air Force Base, Edwards CA, USA}
\ead{ramakanth.munipalli@us.af.mil}

\author[label2]{Karen E.~Willcox}
\ead{kwillcox@oden.utexas.edu}

\begin{abstract}
High-performance computing (HPC) has revolutionized our ability to perform detailed simulations of complex real-world processes.
A prominent contemporary example is from aerospace propulsion, where HPC is used for rotating detonation rocket engine (RDRE) simulations in support of the design of next-generation rocket engines; however, these simulations take millions of core hours even on powerful supercomputers, which makes them impractical for engineering tasks like design exploration and risk assessment.
Data-driven reduced-order models (ROMs) aim to address this limitation by constructing computationally cheap yet sufficiently accurate approximations that serve as surrogates for the high-fidelity model.
This paper contributes a distributed memory algorithm that achieves fast and scalable construction of predictive physics-based ROMs trained from sparse datasets of extremely large state dimension.
The algorithm learns structured physics-based ROMs that approximate the dynamical systems underlying those datasets.
This enables model reduction for problems at a scale and complexity that exceeds the capabilities of standard, serial approaches.
We demonstrate our algorithm's scalability using up to $2,048$ cores on the Frontera supercomputer at the Texas Advanced Computing Center.
We focus on a real-world three-dimensional RDRE for which one millisecond of simulated physical time requires one million core hours on a supercomputer.
Using a training dataset of $2,536$ snapshots each of state dimension $76$ million, our distributed algorithm enables the construction of a predictive data-driven reduced model in just $13$ seconds on $2,048$ cores on Frontera.
\end{abstract}


\end{frontmatter}

\begin{tikzpicture}[remember picture,overlay]
    \node[anchor=base west, scale=0.92] at (current page footer area.south west) {Distribution Statement A: Approved for Public Release; Distribution is Unlimited. PA$\#$ AFRL-2024-1411};
\end{tikzpicture}

\thispagestyle{empty} 

\section{Introduction}
Scientists and engineers leverage data-driven reduced-order models (ROMs) for efficient predictions and decision-making across a wide range of applications underpinned by complex physical phenomena, such as design of next-generation propulsion devices, assessing the impact of turbulent transport in fusion devices, and real-time control of wind farms. 
The advent of petascale computing and, more recently, exascale-capable machines is revolutionizing our ability to conduct numerical simulations of complex real-world problems~\cite{At23}.
This has opened doors to previously unimaginable levels of realism, enabling predictive simulations with billions of degrees of freedom in fusion plasmas~\cite{BCS23, bird2021vpic} or combustion processes~\cite{Ba20,PeleSoftware}, yet the dynamics in such simulations occur across multiple spatiotemporal scales and require substantial computational resources that can often exceed millions of core hours on supercomputers.
The role of ROMs in enabling rapid yet accurate prediction and simulation-based design thus remains as important as ever; however, scalable methods to construct ROMs have not kept pace with the increased scale and resolution made possible by high-performance computing (HPC) in such applications.  

Driven by increased data availability and computing power, various data-driven methods have emerged for scientific applications, including inferring models from data~\cite{BPK16, PW16, schaeffer2017learning}, data-driven discretizations for partial differential equations (PDEs)~\cite{Ba19}, and physics-informed machine learning methods~\cite{Ka21}, to name only a few.
Moreover, a number of recent efforts were dedicated toward data-driven surrogate modeling~\cite{Az24} and model reduction of nonlinear systems~\cite{Ax23, Be20, Ce22, GWW23, GPD17, GB24, Hi20, Pap22, Pe20}, including Operator Inference (OpInf)~\cite{KPW24, PW16}.
OpInf incorporates the structure of the governing equations into the learning problem and constructs structure-preserving ROMs for systems with polynomial structure from data.
However, standard, serial ROM methods are unable to efficiently handle the large-scale training datasets generated by modern HPC simulations, where the state dimension numbers in the many millions (or even billions) and the training dataset---although sparse in time---consumes terabytes of memory.
There is thus a need for novel ROM approaches based on the principles of modern computational science, namely scalable algorithms rooted in physical principles capable of running efficiently on powerful computing systems.
The distributed memory (i.e., message-passing based) ROM algorithm proposed in this work aims to fulfill this need.

Algorithms and software tools for model reduction in large-scale applications have advanced considerably in recent years.
A prominent class of model reduction methods is based on proper orthogonal decomposition (POD)~\cite{BHL93}.
Given the large computational cost of computing POD bases using large datasets, several works have focused on parallel POD approximations via either the singular value decomposition (SVD)~\cite{GvL96} or the method of snapshots~\cite{Si87}.
Beattie et al.~\cite{Be06} proposed an approximate iterative parallel algorithm and Wang, McBee, and Iliescu~\cite{WMI16} formulated an approximate distributed-memory approach based on the method of snapshots for approximating the POD basis.
The recent work~\cite{Ke23} formulated a high-performance solver for computing partial SVD spectra.
The randomized SVD~\cite{HMT11} provides an alternative approach to approximating the POD basis for large datasets.
Reference~\cite{Sayadi2016} employed the parallel QR-decomposition algorithm for tall-and-skinny data matrices proposed in~\cite{De12} to parallelize Dynamic Mode Decomposition (DMD)~\cite{Ku16, Sc10, Tu14}.
Streaming approaches offer a complementary perspective by avoiding the need to fully load the data matrix into memory.
Levy and Lindenbaum~\cite{LL98} proposed an incremental procedure for computing the POD from streaming data.
This procedure was leveraged in Ref.~\cite{schwerdtner2024onlinelearningquadraticmanifolds}, which formulated a streaming, greedy approach for nonlinear dimensionality reduction based on quadratic manifolds.
Nonetheless, when the snapshot dimension is extremely large, performing a streaming approach without also using distributed memory computing will likely be infeasible on a single computer due to the significant memory requirements of loading and processing the data streams.
Complementing methodological improvements, recent efforts have also developed software tools to facilitate efficient model reduction in large-scale applications.
Examples include \texttt{libROM}~\cite{libROM}, a \texttt{C++} library for data-driven approaches including POD, DMD, and projection-based and hyper-reduced intrusive ROMs, \texttt{Pressio}~\cite{rizzi2021pressioenablingprojectionbasedmodel}, an open-source project aimed at providing intrusive model reduction capabilities to large-scale application codes, \texttt{PyMOR}~\cite{MRS16}, a \texttt{Python} library for model order reduction algorithms such as POD, DMD, reduced basis methods, and system-theoretic methods for linear time-invariant systems, which also supports distributed memory parallelization, and the \texttt{PySPOD} library~\cite{Ro24} for parallel spectral POD.

In this paper, we propose a distributed OpInf (dOpInf) algorithm that incorporates HPC into the data-driven learning process to enable rapid and scalable learning of structured, physics-based ROMs for problems at a scale and complexity that exceed the capabilities of standard serial processing approaches for model reduction.
In contrast to existing parallel intrusive methods that require some level of access to the high-fidelity code or approaches that parallelize only specific stages of the reduced modeling process (such as approximating the POD basis), our formulation represents a fully distributed memory workflow for nonlinear, physics-based model reduction of complex applications.
The proposed model reduction workflow enables (1)~efficient data transformations and dimensionality reduction of large datasets with extremely large state dimension without explicitly having to compute the POD basis and without introducing approximations, and (2)~the non-intrusive learning of structured, physics-based ROMs that approximate the dynamics underlying those datasets.
We note that the elements of our distributed algorithm are transferable to other data-driven reduced modeling approaches such as DMD and quadratic manifolds~\cite{barnett2022quadratic, GWW23}, and to parametric ROMs that embed parametric dependence using the strategies surveyed in~\cite{BGW15}.
The scalability and prediction capabilities of the proposed dOpInf algorithm are assessed in a real-world three-dimensional unsteady rotating detonation rocket engine (RDRE) scenario with $76$ million degrees of freedom. 
The training dataset amounts to $1.4$ TB, which prohibits using standard serial methods for model reduction.
In addition to contributing the dOpInf algorithm, our work also represents a novel contribution to model reduction for RDREs, addressing unprecedented levels of complexity in terms of data size and physical coupling.
We demonstrate the scalability of our algorithm using up to $2,048$ cores on the Frontera supercomputer at the Texas Advanced Computing Center (TACC)~\cite{St20}.
On this system, our method constructs a predictive physics-based ROM from the large dataset in just $13$ seconds.
This is in contrast to many existing data-driven reduced and surrogate modeling approaches that can require significant processing times of large datasets and a large snapshot dimension.
The resulting ROM is $90,000$ times faster to evaluate than the original high-fidelity simulation, paving the way toward exploring ROM-based design and quantification of uncertainty.
An implementation, including the detailed tutorial from Ref.~\cite{Fa25}, can be found at \url{https://github.com/ionutfarcas/distributed_Operator_Inference}.

The remainder of this paper is organized as follows.
Section~\ref{sec:setup} summarizes the general setup for high-fidelity nonlinear simulations.
Section~\ref{sec:dOpInf} details the proposed dOpInf algorithm that enables fast and scalable learning of physics-based ROMs for complex applications with training datasets with extremely large state dimension.
We present the scalability results and prediction capabilities of the proposed distributed algorithm in a large-scale RDRE scenario in Sec.~\ref{sec:results}.
Section~\ref{sec:conclusions} concludes the paper.

\section{Setup for high-fidelity nonlinear physics-based simulations} \label{sec:setup}

Consider a complex physical process whose temporal dynamics over $[t_{\textrm{init}}, t_{\textrm{final}}]$ are described by the high-dimensional dynamical system 
\begin{equation} \label{eq:FOM_general}
    \dot{\mathbf{s}} = \mathbf{f}(t, \mathbf{s}), \quad \mathbf{s}(t_{\textrm{init}})=\mathbf{s}_{\mathrm{init}}.
\end{equation}
Here, $\mathbf{s}(t) \in \mathbb{R}^{n}$ denotes the spatially discretized vector of physical state variables at time $t \in [t_{\textrm{init}}, t_{\textrm{final}}]$, $n = n_x n_s \in \mathbb{N}$ is the dimension of the discretized state space, where $n_x \in \mathbb{N}$ denotes the number of degrees of freedom used to discretize the underlying physical domain and $n_s \in \mathbb{N}$ denotes the number of physical state variables (e.g., pressure, velocity, chemical concentrations etc.~in RDRE simulations), $\mathbf{s}_{\mathrm{init}}$ is a specified initial condition, and $\mathbf{f} : [t_{\textrm{init}}, t_{\textrm{final}}] \times \mathbb{R}^{n} \rightarrow \mathbb{R}^{n}$ is a nonlinear function that defines the time evolution of $\mathbf{s}$.

Our starting point is a given set of training data that represent observations of the state of the system, where the training data might be generated by experiments or by simulation, or by a combination of the two.
In many cases, and in the real-world RDRE scenario considered in this paper, the data come from a simulation that solves the physical equations governing the spatiotemporal evolution of the system state. 
These governing equations take the form of PDEs encoding physical conservation laws and constitutive relationships. 
When discretized in space, these equations are written as an $n$-dimensional system of ordinary differential equations (ODEs) as in~\eqref{eq:FOM_general}, representing the system of discretized PDEs over the computational domain denoted by $\Omega \subset \mathbb{R}^d$, where typically $d = 2, 3$.
In such cases, $n$ is large (e.g., in $O(10^6) - O(10^9)$), because it scales with the dimension of the PDE spatial discretization, $n_x$.

We focus on ROMs utilized for predictions over a time horizon $[t_{\textrm{init}}, t_{\textrm{final}}]$, with $t_{\textrm{init}}$ denoting the initial time and $t_{\textrm{final}}$ denoting the final time.
The training dataset comprises $n_t \in \mathbb{N}$ large-scale state vectors or snapshots at $n_t$ time instants over a training horizon $[t_{\textrm{init}}, t_{\textrm{train}}]$ with $t_{\textrm{train}} < t_{\textrm{final}}$.
The state solution at time instant $t_k \in [t_{\textrm{init}}, t_{\textrm{train}}]$ is referred to as the $k$th snapshot and is denoted by $\mathbf{s}_k$.
The $n_t$ snapshots are collected into a large-scale snapshot matrix $\mathbf{S} \in \mathbb{R}^{n \times n_t}$ with $\mathbf{s}_k$ as its $k$th column:
\begin{equation*}
    \mathbf{S} =
     \begin{bmatrix}
\vert & \vert & & \vert & & \vert\\
     \mathbf{s}_1 &
     \mathbf{s}_2 &
     \ldots &
     \mathbf{s}_k &
     \ldots &
     \mathbf{s}_{n_t}\\
     \vert & \vert & & \vert & & \vert
     \end{bmatrix}.
\end{equation*}
In large-scale applications, $n_x$ (and therefore $n = n_s n_x$) is orders of magnitude larger than $n_t$.
In other words, the row dimension of $\mathbf{S}$ is significantly larger than its column dimension.
Given $\mathbf{S}$ and the governing equations~\eqref{eq:FOM_general}, our goal is to efficiently and scalably learn a predictive data-driven ROM over the target time horizon $[t_{\textrm{init}}, t_{\textrm{final}}]$.
However, doing so via standard serial methods would be computationally expensive and limited by memory storage and bandwidth.
Our proposed distributed workflow, in contrast, enables a fast and scalable workflow for constructing physics-based ROMs on distributed memory machines.

\section{dOpInf: A new distributed computing algorithm for fast and scalable learning of nonlinear physics-based reduced models} \label{sec:dOpInf}
This section presents the proposed distributed Operator Inference (dOpInf) algorithm.
dOpInf represents a distributed memory workflow for physics-based model reduction of complex applications with extremely large state dimension.
Section~\ref{subsec:dOpInf_preproc} presents the first step in dOpInf focused on distributed data manipulations and dimensionality reduction, followed by the distributed learning of the reduced model operators in Sec.~\ref{subsec:dOpInf_learning}, and the distributed postprocessing of the reduced solution in Sec.~\ref{subsec:dOpInf_postprocessing}.
We provide a summary in Sec.~\ref{subsec:dOpInf_summary}.

\subsection{Distributed computation of data transformations and dimensionality reduction} \label{subsec:dOpInf_preproc}
The OpInf methodology requires the generation of a low-dimensional representation of the snapshot data. 
These low-dimensional data are then used to infer the reduced model operators that define the physics-based ROM. 
Generating the low-dimensional snapshot representation entails two key steps: data transformations and dimensionality reduction. 
Data transformations are used to improve numerical conditioning through centering and scaling, and to expose polynomial structure in nonlinear problems \cite{Qi20}.
Dimensionality reduction represents the (transformed) snapshot data in a subspace of reduced dimension $r \ll n$. 
The proposed dOpInf approach introduces a formulation that ensures a scalable distributed computation for each of these steps. 
In particular, by starting from the POD method of snapshots introduced in the seminal work~\cite{Si87} but now typically replaced in POD implementations by the SVD, we formulate a more scalable algorithm for achieving dimensionality reduction without necessitating the computation of the POD basis.

Algorithm~\ref{alg:dOpInf_preprocessing} details our distributed memory formulation of data transformations and dimensionality reduction.
The inputs are the number of compute cores $p \in \mathbb{N}_{\geq 2}$, the high-dimensional snapshot dimension $n$, the number of training snapshots $n_t$, the reduced dimension of the target ROM $r \in \mathbb{N}$, and a variable \emph{transform} that indicates whether data transformations are employed.
In the following, we denote by $i=1, 2,\ldots,p$ the index of the $i$th core.
\begin{algorithm}
\caption{dOpInf data transformations and dimensionality reduction on $p$ compute cores}\label{alg:dOpInf_preprocessing}
\textbf{Input}: $p \geq 2$, $n$, $n_t$, $r$, \emph{transform} \\
\textbf{Output}: projected transformed snapshots $\hat{\mathbf{Q}} \in \mathbb{R}^{r \times n_t}$ 
\begin{algorithmic}[1]
\State \parbox[t]{0.90\textwidth}{\emph{divide} the snapshot dimension $n$ into $p$ parts $n_1, n_2, \ldots, n_p$ such that $\sum_{i=1}^p n_i = n$} 
\For{$i \gets 1$ to $p$ in parallel\strut}
\State \emph{load} snapshot components $\mathbf{S}_i \in \mathbb{R}^{n_i \times n_t}$ from disk
\If{\emph{transform}}
    \State \parbox[t]{0.90\textwidth}{apply \emph{data transformations} to $\mathbf{S}_i$ such as lifting, centering, and scaling, to obtain $\mathbf{Q}_i \in \mathbb{R}^{m_i \times n_t}$ such that $\sum_{j=1}^p m_j = m \geq n$\strut}
\Else
    \State continue with $\mathbf{Q}_i = \mathbf{S}_i$  and set $m_i = n_i$
\EndIf
\State \emph{compute} $\mathbf{D}_i = \mathbf{Q}_i^\top \mathbf{Q}_i \in \mathbb{R}^{n_t \times n_t}$
\State \parbox[t]{0.90\textwidth}{\emph{compute} and \emph{broadcast} $\mathbf{D} = \sum_{j=1}^p \mathbf{D}_j$ via a collective reduction \strut}

\State \parbox[t]{0.90\textwidth}{\emph{compute} the first $r$ eigenpairs $\{(\lambda_k, \mathbf{u}_k)\}_{k=1}^{r}$ of $\mathbf{D}$ and arrange them s.~t.~$\lambda_1 \geq \lambda_2 \geq \ldots \geq \lambda_{r}$\strut}

\State \parbox[t]{0.90\textwidth}{\emph{compute} the corresponding global POD singular values $\sigma_k = \sqrt{\lambda_k}$ for $k = 1, 2, \ldots, r$\strut}

\State \parbox[t]{0.90\textwidth}{\emph{compute} $\mathbf{T}_r = \mathbf{U}_r \bm{\Lambda}_r^{-\frac{1}{2}} \in \mathbb{R}^{n_t \times r}$, where $\mathbf{U}_r = \begin{bmatrix} \mathbf{u}_1 \vert \mathbf{u}_2 \vert \dots \vert \mathbf{u}_r \end{bmatrix}$ and $\bm{\Lambda}_r = \mathrm{diag}(\lambda_1, \lambda_2, \ldots, \lambda_r)$\strut}

\State \parbox[t]{0.90\textwidth}{\emph{compute} the global representation of the transformed data in the low-dimensional space spanned by the POD basis vectors as $\hat{\mathbf{Q}} = \mathbf{T}_r^\top \mathbf{D} \in \mathbb{R}^{r \times n_t}$}\strut

\EndFor
\end{algorithmic}
\end{algorithm}

The first step is to distribute the snapshot data into $p$ non-overlapping row blocks $\mathbf{S}_i \in \mathbb{R}^{n_i \times n_t}$ with one block per core such that $\sum_{i = 1}^p n_i = n$.
The splitting strategy is user-defined and depends on factors such as the target parallel architecture and necessary data transformations, as we will discuss below.
A straightforward approach suitable for homogeneous architectures is to use blocks of size $n_i = \lfloor n/p \rfloor$; when $p$ does not divide $n$ exactly, the remaining $n - \lfloor n/p \rfloor$ rows can be further distributed among the $p$ compute cores.
Loading $\mathbf{S}_i$ into the memory of the $i$th core requires $\mathcal{O}(n_i n_t)$ bytes.
\begin{remark} \label{remark:data_format}
    To ensure scalable access to the training dataset across $p$ compute cores, large datasets should be stored in a format like HDF5 or NetCDF, which supports efficient parallel I/O operations, and on high throughput partitions such as the \texttt{scratch} partition on a supercomputer.
    A comprehensive understanding of the underlying file system is also crucial for maximizing performance.
    We note that saving a large-scale dataset in a single file can hinder scalability when an increasing number of cores attempt to access it simultaneously. 
    This can be mitigated by file partitioning (i.e., the dataset is divided into multiple files, allowing scalable parallel reading operations across different compute cores) or distributed reading and broadcasting (i.e., one core reads the data and broadcasts it to all other cores), keeping in mind that this strategy might be more time consuming than file partitioning for very large datasets.
\end{remark}

The application of data transformations is determined by the input \emph{transform}. 
While often glossed over in the literature, most complex applications of dimensionality reduction via POD require centering and scaling of the snapshot data.
The importance of centering and scaling for obtaining accurate OpInf ROMs is particularly observed in problems with multiple physical state variables (e.g., pressure, velocity, chemical concentrations in RDRE simulations), each with differing physical scales~\cite{Fa23, QFW21, Sw20}.
One strategy, which we also employ in the RDRE scenario under consideration, is to center each discretized state variable in $\mathbf{S}_i$ by its mean over the training time horizon and then scale it by its maximum absolute value to ensure that the scaled variables do not exceed $[-1, 1]$.
To efficiently compute centered and scaled snapshots, we distribute the full snapshot dataset over the full computational domain $\Omega$ into $p$ non-overlapping subdomains $\Omega_1, \Omega_2, \ldots, \Omega_p \in \mathbb{R}^d$ satisfying $\Omega = \cup_{i = 1}^p \Omega_i$ such that each compute core gets all discrete state variables for a subdomain $\Omega_i$.
Note that this splitting strategy aligns with traditional non-overlapping domain decomposition methods.
This allows for independent centering calculations on each core, eliminating communication needs. 
Scaling parameters, which are typically global across the training horizon, can be computed locally on each core followed by an inexpensive collective communication step to compute the global results.
In general, the exact computational costs of centering and scaling are transformation dependent, but they are usually low compared to the remaining preprocessing costs.
In addition, they can be performed in-place on each $\mathbf{S}_i$, requiring no additional memory.
The OpInf approach might also employ lifting and other variable transformations~\cite{Qi20}. 
By exploiting the knowledge of a system’s governing equations, these transformations seek to find a new coordinate system where the system dynamics exhibit a polynomial structure.
One example used in computational fluid dynamics involves expressing the governing equations in terms of specific volume variables instead of the conventional conservative variables, resulting in a quadratic structure in the transformed governing equations.
In other examples, the transformations might involve augmenting the system's physical state with additional auxiliary variables, thus lifting the physical variables to a higher dimensional coordinate system.
Lifting and other variable transformations are applied entry-wise in each $\mathbf{S}_i$ and do not generally require communication.

If data transformations are necessary, the transformed snapshot matrices for each compute core are denoted as $\mathbf{Q}_i \in \mathbb{R}^{m_i \times n_t}$, where $m_i$ denotes the dimensions of the transformed $i$th snapshot partition, such that $\sum_{i = 1}^p m_i = m \geq n$; $m$ exceeds $n$ when employing lifting transformations that introduce auxiliary state variables. 
In this case, the number of lifted (or transformed) state variables, $m_s \in \mathbb{N}$, exceeds the original number of physical state variables, $n_s$.
The original partitions $\mathbf{S}_i$ can be discarded from memory, and $\mathbf{Q}_i$ are used for the ensuing calculations.
If data transformations are not required, we set $\mathbf{Q}_i=\mathbf{S}_i$ and $m_i = n_i$. 
Under these circumstances, any snapshot partitioning scheme suffices provided that the entire dataset is non-overlappingly distributed among the $p$ cores.

We next compute the representation of the transformed snapshots in the low-dimensional subspace spanned by the rank-$r$ POD basis vectors.
We start from the method of snapshots~\cite{Si87}, as this provides several computational and memory advantages in our context.
We note that the sequential version of the method of snapshots offers, in principle, a path to handle large datasets by processing subsets of at least two snapshots at a time.
However, its sequential nature and potential for redundancy (e.g., reloading the full dataset for scaling) limit its efficiency.

We compute the Gram matrices $\mathbf{D}_i = \mathbf{Q}_i^\top \mathbf{Q}_i \in \mathbb{R}^{n_t \times n_t}$ on each compute core, which involves a dense matrix-matrix multiplication requiring $\mathcal{O}(m_i n_t^2)$ operations and $\mathcal{O}(n_t^2)$ bytes of main memory.
We then compute their summation $\mathbf{D} = \sum_{i=1}^p \mathbf{D}_i$ via a collective parallel reduction that ensures that all cores have $\mathbf{D}$ in their local memory, which requires $\mathcal{O}(n_t^2)$ computations and $\mathcal{O}(n_t^2)$ additional bytes of memory on each core. 
Since the full snapshot data was distributed non-overlappingly, the following lemma shows that $\mathbf{D}$ is equal to the global Gram matrix $\mathbf{Q}^\top \mathbf{Q}$. 
\begin{lemma} \label{lemma:MoS_split}
    Let $\mathbf{A} \in \mathbb{R}^{q \times k}$ be a matrix and let $\mathbf{A}_i \in \mathbb{R}^{q_i \times k}$ with $i = 1, 2, \ldots, p$ be $p$ non-overlapping row blocks such that $\sum_{i=1}^p q_i = q$.
    Then, $\mathbf{A}^\top \mathbf{A} = \sum_{i=1}^p \mathbf{A}_i^\top \mathbf{A}_i$.
\begin{proof}
    Let $\mathbf{P}_i \in \mathbb{R}^{q \times q_i}$ be prolongation matrices with entries $p_{i, j} \in \{0, 1\}$ such that $\mathbf{A} = \sum_{i=1}^p \mathbf{P}_i \mathbf{A}_i$. 
    The conclusion follows from the fact that $\mathbf{P}_i^\top \mathbf{P}_i = \mathbf{I}_{q_i \times q_i}$ and, since the row blocks are non-overlapping,  $\mathbf{P}_j^\top \mathbf{P}_i = \mathbf{0}_{q_j \times q_i}, \quad \forall j \neq i$.
\end{proof}
\end{lemma} 

Each core proceeds by computing the partial eigendecomposition $\{(\lambda_k, \mathbf{u}_k)\}_{k=1}^{r}$ of the symmetric positive semi-definite matrix $\mathbf{D}$, where $\lambda_k$ are the real and non-negative eigenvalues, and $\mathbf{u}_k \in \mathbb{R}^{n_t}$ denote the corresponding eigenvectors.
The computational cost of this decomposition depends on the value of $r$ and the employed approach.
If $r$ is small, efficient methods such as those based on the Lanczos iterative algorithm~\cite{GvL96} can be employed.
The $r$ eigenpairs necessitate $\mathcal{O}(r n_t)$ bytes of memory on each core.
We must ensure that they are arranged such that $\lambda_1 \geq \lambda_2 \geq \ldots \geq \lambda_{r}$.

Consider the thin SVD~\cite{GvL96} of $\mathbf{Q} = \mathbf{V} \bm{\Sigma} \mathbf{W}^\top$.
$\mathbf{V} \in \mathbb{R}^{m \times n_t}$ contains the left singular vectors, $\bm{\Sigma} \in \mathbb{R}^{n_t \times n_t}$ is a diagonal matrix containing the singular values of $\mathbf{Q}$ in non-decreasing order $\sigma_1 \geq \sigma_2 \geq \ldots \geq \sigma_{n_t}$, where $\sigma_j$ denotes the $j$th singular value, and $\mathbf{W} \in \mathbb{R}^{n_t \times n_t}$ contains the right singular vectors.
$\mathbf{V}$ and $\mathbf{W}$ are semi-orthogonal, that is, $\mathbf{V}^\top \mathbf{V} = \mathbf{W}^\top \mathbf{W} = \mathbf{I}_{n_t}$, where $\mathbf{I}_{n_t} \in \mathbb{R}^{n_t \times n_t}$ denotes the identity matrix.
Since 
\begin{equation*}
    \mathbf{D} = \mathbf{Q}^\top \mathbf{Q} = \mathbf{W} \mathbf{\Sigma} \mathbf{V}^\top \mathbf{V} \mathbf{\Sigma} \mathbf{W}^\top = \mathbf{W} \mathbf{\Sigma}^2 \mathbf{W}^\top\;,
\end{equation*}
it follows that $\mathbf{D} \mathbf{W} = \mathbf{W} \mathbf{\Sigma}^2$, where we used the semi-orthogonality of $\mathbf{W}$.
This implies that the eigenvalues of $\mathbf{D}$ are the squared singular values of $\mathbf{Q}$ and the eigenvectors of $\mathbf{D}$ are equivalent (in terms of spanned subspaces) to the right singular vectors of $\mathbf{Q}$.
Hence, the first $r$ global POD singular values can be computed as $\sigma_k = \sqrt{\lambda_k}, \quad k = 1, 2, \ldots, r$, which requires $\mathcal{O}(r)$ operations.
Furthermore, from the thin SVD of $\mathbf{Q}$ we have that $\mathbf{V} = \mathbf{Q} \mathbf{W} \bm{\Sigma}^{-1}$, which means that the left singular vectors in $\mathbf{V}$ can be equivalently expressed as
\begin{equation} \label{eq:MoS_SVD_left_singular_vectors}
    \mathbf{V} = \mathbf{Q} \mathbf{U} \bm{\Lambda}^{-\frac{1}{2}}, 
\end{equation}
where $\mathbf{U} = \begin{bmatrix} \mathbf{u}_1 \vert \mathbf{u}_2 \vert \dots \vert \mathbf{u}_{n_t} \end{bmatrix}$ and $\bm{\Lambda}_r = \mathrm{diag}(\lambda_1, \lambda_2, \ldots, \lambda_{n_t})$.

The user can also choose $r$ during the execution of Algorithm~\ref{alg:dOpInf_preprocessing} based on an energy criterion such as
\begin{equation} \label{eq:ret_energy}
    \frac{\sum_{k=1}^{r} \sigma_k^2}{\sum_{k=1}^{n_t} \sigma_k^2} \geq E_{\%},
\end{equation}
where $E_{\%}$ is a user-specified energy truncation threshold (e.g., $E_{\%} = 95\%$ or $E_{\%} = 99\%$).
In this case, we must compute all eigenpairs of $\mathbf{D}$, arrange them such that the eigenvalues are in increasing order, and choose $r$ based on~\eqref{eq:ret_energy}.
The first $r$ eigenpairs are then used in the next steps in Algorithm~\ref{alg:dOpInf_preprocessing}.  
Computing all eigenpairs of a symmetric positive definite matrix requires $\mathcal{O}(n_t^3)$ operations and $\mathcal{O}(n_t^2)$ bytes of memory for the $n_t$ eigenpairs.

In the next step (step $12$), each core computes $\mathbf{T}_r = \mathbf{U}_r \bm{\Lambda}_r^{-\frac{1}{2}} \in \mathbb{R}^{n_t \times r}$, where $\mathbf{U}_r = \begin{bmatrix} \mathbf{u}_1 \vert \mathbf{u}_2 \vert \dots \vert \mathbf{u}_r \end{bmatrix}$ and $\bm{\Lambda}_r = \mathrm{diag}(\lambda_1, \lambda_2, \ldots, \lambda_r)$, which requires $\mathcal{O}(r n_t)$ computations and $\mathcal{O}(r n_t)$ bytes of memory for the result.
At this point we have all the ingredients to compute the global representation of the transformed snapshots in the low-dimensional subspace spanned by the rank-$r$ POD basis vectors without explicitly requiring the basis.
The standard, sequential approach to compute the low-dimensional representation $\hat{\mathbf{Q}} \in \mathbb{R}^{r \times n_t}$ of the high-dimensional transformed data in OpInf is
\begin{equation}\label{eq:std_OpInf_data_projection}
    \hat{\mathbf{Q}} = \mathbf{V}_r^\top \mathbf{Q},
\end{equation}
where $\mathbf{V}_r \in \mathbb{R}^{m \times r}$ denotes the global rank-$r$ POD basis.
In the standard OpInf formulation~\cite{KPW24,PW16}, $\mathbf{V}_r$ is computed from the thin SVD of the full snapshot matrix $\mathbf{Q}$ by taking the first $r \ll m$ columns of $\mathbf{V}$, that is, the left singular vectors corresponding to the $r$ largest singular values.
While~\eqref{eq:std_OpInf_data_projection} is easily parallelizable, it would involve several steps: (i) compute the components $\mathbf{V}_{r, i} \in \mathbb{R}^{m_i \times r}$ of the POD basis on each core (which can be done by leveraging~\eqref{eq:MoS_SVD_left_singular_vectors} or the thin SVD of $\mathbf{Q}_i$, for example), (ii) compute $\hat{\mathbf{Q}}_i = \mathbf{V}_{r, i}^\top \mathbf{Q}_i$ on each core, which would require $\mathcal{O}(m_i n_t r)$ operations and $\mathcal{O}(r n_t)$ bytes of memory, and (iii) compute $\hat{\mathbf{Q}}$ via a parallel reduction in terms of $\hat{\mathbf{Q}}_1, \hat{\mathbf{Q}}_2, \ldots, \hat{\mathbf{Q}}_p$, which would require communication.
In contrast, from Eq.~\eqref{eq:MoS_SVD_left_singular_vectors}, the global rank-$r$ POD basis can be equivalently computed as 
\begin{equation} \label{eq:MoS_POD_basis}
    \mathbf{V}_r = \mathbf{Q} \mathbf{W}_r \bm{\Sigma}_r^{-1} = \mathbf{Q} \mathbf{U}_r \bm{\Lambda}_r^{-\frac{1}{2}}. 
\end{equation}
This, in turn means that 
\begin{equation} \label{eq:MoS_OpInf_projection}
    \hat{\mathbf{Q}} = \left( \mathbf{Q} \mathbf{U}_r \bm{\Lambda}_r^{-\frac{1}{2}}\right)^\top  \! \! \! \mathbf{Q} = \mathbf{T}_r^\top \mathbf{Q}^\top \mathbf{Q} = \mathbf{T}_r^\top \mathbf{D}.
\end{equation}
Therefore, by starting from the POD method of snapshots, $\hat{\mathbf{Q}}$ can be computed in terms of the small matrices $\mathbf{T}_r$ and $\mathbf{D}$ at a much lower cost of $\mathcal{O}(r n_t^2)$ operations and $\mathcal{O}(r n_t)$ bytes of local main memory without necessitating the computation of the POD basis.

The total cost of Algorithm~\ref{alg:dOpInf_preprocessing} per core is in $\mathcal{O}(m_i n_t^2 + n_t^2 + rn_t + rn_t^2) = \mathcal{O}((m_i + r + 1) n_t^2) = \mathcal{O}(m_i n_t^2)$ since $r$ is generally smaller than $m_i$, plus the cost of computing the eigendecomposition of $\mathbf{D}$ and the global POD singular values from the eigenvalues of $\mathbf{D}$.
The corresponding memory requirements are in $\mathcal{O}(m_i n_t + 2n_t^2 + 2rn_t)$ bytes plus the memory requirements for holding the eigenpairs of $\mathbf{D}$.

\begin{remark}\label{remark:collective_reduction_Algo_1}
    Step $9$ in Algorithm~\ref{alg:dOpInf_preprocessing} computes and broadcasts $\mathbf{D}$ to all other cores via a collective reduction, and subsequently all cores compute the reduced eigendecomposition of $\mathbf{D}$ as well as the matrices $\mathbf{T}_r$ and $\hat{\mathbf{Q}}$ (steps $10$--$13$).
    This is because the reduced representation of the transformed snapshot data, $\hat{\mathbf{Q}}$, is required by all cores for inferring the reduced model operators via OpInf in the next step, and computing $\hat{\mathbf{Q}}$ depends on $\mathbf{T}_r$ and $\mathbf{D}$ via~\eqref{eq:MoS_OpInf_projection}.
    In addition, as we will discuss in Section~\ref{subsec:dOpInf_postprocessing}, $\mathbf{T}_r$ is generally also needed by all cores for postprocessing the reduced solution produced by the dOpInf ROM.
    These computations are cheap when $n_t$ and $r$ are significantly smaller than $m$, in which case they do not affect the scalability of the algorithm (as demonstrated by our scalability results in Section~\ref{subsec:results_scalability}).
    An alternative approach would be to compute $\mathbf{D}$ at step $9$ using a standard reduction to one core, perform all subsequent computations on that core and broadcast the results to all other compute cores. 
    However, this strategy requires communication and can create a bottleneck, as all remaining computations must wait for the results, which could reduce the overall efficiency of the data dimensionality reduction procedure.
\end{remark}
\begin{remark}\label{remark:scalability_Algo_1}
    Algorithm~\ref{alg:dOpInf_preprocessing} is distributed with respect to the large transformed snapshot dimension, $m$.
    All other steps (steps $10$--$13$) are sequential in the current formulation. 
    When $n_t$ is large, the cost of these steps per core will become dominant, which can impact the scalability of Algorithm~\ref{alg:dOpInf_preprocessing}.
    One solution is to compute $\mathbf{T}_r$ (step $12$) and $\hat{\mathbf{Q}}$ (step $13$) in parallel on the $p$ cores with respect to $n_t$.
    In addition, the full or reduced eigendecomposition of $\mathbf{D}$ at step $10$ can be performed in parallel using high-performance linear algebra libraries like ELPA\footnote{\url{https://elpa.mpcdf.mpg.de/}} (Eigenvalue Solvers for Petaflop Applications)~\cite{Au11}.
    Another solution is to leverage streaming methods: instead of processing all $n_t$ snapshots at once, we would process them in streams of size $n_\ell$ such that $n_{\ell} < n_t$~\cite{LL98}. 
    Future work will explore incorporating streaming approaches into Algorithm~\ref{alg:dOpInf_preprocessing}.
\end{remark}

\subsection{Distributed learning of the reduced model operators} \label{subsec:dOpInf_learning}
We next use the reduced data matrix $\hat{\mathbf{Q}}$ to learn the reduced operators that specify the reduced model.
For example, for a quadratic reduced model
\begin{equation} \label{eq:ROM_quad_time_cont}
    \dot{\hat{\mathbf{q}}} = \hat{\mathbf{A}}\hat{\mathbf{q}} + \hat{\mathbf{H}}\left(\hat{\mathbf{q}} \otimes \hat{\mathbf{q}} \right)  + \hat{\mathbf{c}},
\end{equation}
we must determine the constant, linear, and quadratic reduced operators $\hat{\mathbf{c}} \in \mathbb{R}^{r}, \hat{\mathbf{A}} \in \mathbb{R}^{r\times r}$, and $\hat{\mathbf{H}} \in \mathbb{R}^{r\times r^2}$.
OpInf determines the reduced operators that best match the projected snapshot data in a minimum residual sense by solving the linear least-squares minimization
\begin{equation} \label{eq:OpInf_reg}
    \mathop{\mathrm{argmin}}_{\hat{\mathbf{O}}} \left\lVert \hat{\mathbf{D}}\hat{\mathbf{O}}^{\top} - \dot{\hat{\mathbf{Q}}}^\top \right\rVert_F^2 + \beta_{1} \left(\left\lVert\hat{\mathbf{A}}\right\rVert_F^2 + \left\lVert\hat{\mathbf{c}}\right\rVert_F^2\right) + \beta_{2} \left\lVert\hat{\mathbf{H}}\right\rVert_F^2,
\end{equation}
where $\hat{\mathbf{O}} =
\begin{bmatrix}
\hat{\mathbf{A}} \, \vert \, \hat{\mathbf{H}} \, \vert \, \hat{\mathbf{c}}
\end{bmatrix}
\in \mathbb{R}^{r \times (r + r^2 + 1)}$ denotes the unknown operators,  
$\hat{\mathbf{D}} =
\begin{bmatrix}
\hat{\mathbf{Q}}^\top \, \vert \, \hat{\mathbf{Q}}^\top \odot \hat{\mathbf{Q}}^\top \, \vert \, \hat{\mathbf{1}}_{n_t}
\end{bmatrix}
\in \mathbb{R}^{n_t \times (r + r^2 +1)}$ denotes the OpInf data, $\dot{\hat{\mathbf{Q}}} \in \mathbb{R}^{r \times n_t}$ is the time derivative data, and $F$ denotes the Frobenius norm.
When the high-fidelity model does not provide time derivative data, $\dot{\hat{\mathbf{Q}}}$ must be approximated numerically using finite differences, for example. 
Numerical approximations, however, will likely be inaccurate in applications with severely downsampled snapshots, as it was also observed in the RDRE scenario under consideration.
In such cases, we use the time-discrete version of OpInf~\cite{Fa23} and learn the corresponding reduced operators analogously to Eq.~\eqref{eq:OpInf_reg}, where instead of using $\dot{\hat{\mathbf{Q}}}$, we shift $\hat{\mathbf{Q}}$ by one column the right as is done in DMD, which learns linear discrete-time systems.

The learning step that can benefit from a distributed formulation is the search for the optimal regularization hyperparameter pair $(\beta_1^{\mathrm{opt}}, \beta_2^{\mathrm{opt}} ) \in \mathbb{R}_{> 0}^2$.
These parameters are introduced to address overfitting and accommodate model misspecification and other sources of error~\cite{MHW21}.
Following~\cite{MHW21, QFW21}, the regularization hyperparameters can be found through a grid search involving two nested loops.
The regularization pairs are independent of each other, which means that this search is 
embarrassingly parallel.
The optimal hyperparameters minimize the training error under the constraint that the inferred reduced coefficients have bounded growth within a trial time horizon at least as long as the target time horizon~\cite{QFW21}.

We present the steps to learn the reduced operators using a distributed search for the optimal hyperparameters in Algorithm~\ref{alg:dOpInf_learning}.
In addition to the input parameters from Algorithm~\ref{alg:dOpInf_preprocessing}, we have the sets $\mathcal{B}_1$, $\mathcal{B}_2 \subset \mathbb{R}_{> 0}$ comprising the candidate regularization parameters $\beta_1 \in \mathcal{B}_1$ and $\beta_2 \in \mathcal{B}_2$, $t_{\textrm{trial}} \geq t_{\textrm{final}}$ denoting the end of the trial time horizon, and a percentage $\tau \in (0, 1)$ constraining the growth of the inferred reduced coefficients over the trial horizon.
\begin{algorithm}
\caption{dOpInf reduced operator learning}\label{alg:dOpInf_learning}
\textbf{Input}: $p \geq 2$, $n$, $n_t$, $r$, \emph{transform}, $\mathcal{B}_1$, $\mathcal{B}_2$, $t_{\textrm{trial}}$, $\tau \in (0, 1)$ \\
\textbf{Output}: reduced operators $\hat{\mathbf{c}} \in \mathbb{R}^{r}, \hat{\mathbf{A}} \in \mathbb{R}^{r\times r}$, $\hat{\mathbf{H}} \in \mathbb{R}^{r\times r^2}$
\begin{algorithmic}[1]
\State \emph{split} the regularization parameter pairs in $\mathcal{B}_1 \times \mathcal{B}_2$ into $p$ disjoint subsets $\mathcal{B}_1^{i} \times \mathcal{B}_2^{i}$ of equal size

\State \emph{use} the reduced data matrix $\hat{\mathbf{Q}} \in \mathbb{R}^{r \times n_t}$ computed using Algorithm~\ref{alg:dOpInf_preprocessing}

\For{$i \gets 1$ to $p$ in parallel\strut}

\For{$(\beta_1, \beta_2) \in \mathcal{B}_1^{i} \times \mathcal{B}_2^{i}$ \strut}
\State \parbox[t]{0.90\textwidth}{\emph{infer} the reduced operators using $\{\beta_1, \beta_2 \}$}
\State \parbox[t]{0.90\textwidth}{\emph{compute} the reduced solution over the trial time horizon $[t_{\textrm{init}}, t_{\textrm{trial}}]$}
\EndFor

\State \parbox[t]{0.90\textwidth}{perform a parallel reduction to \emph{find} $(\beta_1^{\mathrm{opt}}, \beta_2^{\mathrm{opt}} )$ that minimize the training error and 
ensure that the maximum deviation from the training mean of the inferred reduced coefficients over the trial horizon stays within $\tau$ of the maximum deviation from the mean over training}

\State \parbox[t]{0.90\textwidth}{\emph{determine} the index $i^{\mathrm{opt}}$ of the compute core where $(\beta_1^{\mathrm{opt}}, \beta_2^{\mathrm{opt}} )$ resides}

\If{$i = i^{\mathrm{opt}}$}
\State \parbox[t]{0.90\textwidth}{the target reduced operators $\{\hat{\mathbf{c}}, \hat{\mathbf{A}}, \hat{\mathbf{H}} \}$ are the ones \emph{inferred} using $(\beta_1^{\mathrm{opt}}, \beta_2^{\mathrm{opt}})$}
\EndIf
\EndFor
\end{algorithmic}
\end{algorithm}
For convenience, we choose $\mathcal{B}_1$ and $\mathcal{B}_2$ such that the cardinality $B \in \mathbb{N}$ of $\mathcal{B}_1 \times \mathcal{B}_2$ is divisible by $p$. 
We start by splitting the regularization pairs in $\mathcal{B}_1 \times \mathcal{B}_2$ into $p$ disjoint subsets of equal size $B/p$.
For each regularization pair in $\mathcal{B}_1^{i} \times \mathcal{B}_2^{i}$, we first infer the reduced operators using Eq.~\eqref{eq:OpInf_reg} and then compute the reduced solution over the trial time horizon $[t_{\textrm{init}}, t_{\textrm{trial}}]$.
Since~\eqref{eq:OpInf_reg} depends on the reduced dimension $r$ and the number of training snapshots $n_t$, solving it is computationally cheap in general and can be done utilizing standard least-squares solvers~\cite{GvL96}.
Computing the dOpInf reduced solution is also cheap; in fact, in large-scale applications, this is orders of magnitude faster than the high-fidelity model.
The distributed partitioning of the candidate regularization pairs therefore decreases the cost of the grid search from $B$ least-squares solves~\eqref{eq:OpInf_reg} and reduced solution calculations in the sequential approach to $B/p$ such calculations per core.

The optimal regularization pair $(\beta_1^{\mathrm{opt}}, \beta_2^{\mathrm{opt}})$ is found via a parallel reduction that minimizes the error between the computed reduced solutions and the given projected training data $\hat{\mathbf{Q}}$ (this error can be, for example, the mean-squared error) subject to a constraint. 
The maximum deviation of the inferred reduced coefficients over the trial horizon from the mean of the reduced-order coefficients over training must stay within $\tau$ of the maximum deviation from the mean over training.
We then determine the index $i^{\mathrm{opt}}$ of the compute core where the reduced operators inferred using the optimal regularization pair reside. 

\subsection{Distributed postprocessing of the reduced solution}~\label{subsec:dOpInf_postprocessing}
Finally, we postprocess the reduced solution produced by the dOpInf ROM constructed using the optimal regularization parameter pair $(\beta_1^{\mathrm{opt}}, \beta_2^{\mathrm{opt}})$.
Let $\Tilde{\mathbf{Q}} \in \mathbb{R}^{r \times n_p}$ denote the matrix containing the reduced solution at $n_p \in \mathbb{N}$ target time instants.
For example, $n_p$ can represent user-defined, representative time instants or all time instants over the time horizon of interest $[t_{\textrm{init}}, t_{\textrm{final}}]$.

In many cases, postprocessing involves computing low-dimensional output quantities of interest that are relevant to the problem at hand (e.g., integrated forces, average field values, state values over a specific region of interest, etc.). 
In some cases, postprocessing the reduced solution might involve mapping it back to the original high-dimensional space.
This would be the case, for example, if the user desired visualization of the full state fields, such as pressure or temperature.
In these cases, in order to map from the reduced solution  $\Tilde{\mathbf{Q}}$ to the desired full-state field, we must compute the corresponding components of the rank-$r$ POD basis (noting that our proposed dOpInf algorithm has thus far avoided computing the POD basis).
To achieve this POD basis computation efficiently, we distribute it across cores.
Let $\mathbf{V}_{r, i} \in \mathbb{R}^{m_i \times r}$ denote the components of the global POD basis on the $i$th core, where $i=1, 2, \ldots, p$.
From Eq.~\eqref{eq:MoS_POD_basis}, we can compute $\mathbf{V}_{r, i}$ as
\begin{equation} \label{eq:MoS_POD_basis_comp}
   \mathbf{V}_{r, i} = \mathbf{Q}_i \mathbf{T}_r.
\end{equation}
Equation~\eqref{eq:MoS_POD_basis_comp} requires $\mathcal{O}(m_i n_t r)$ operations and $\mathcal{O}(m_i r)$ bytes of memory for $\mathbf{V}_{r, i}$ on the $i$th core.
Lifting the reduced solution back to the high-dimensional space is done by computing $\mathbf{V}_{r, i} \Tilde{\mathbf{Q}} \in \mathbb{R}^{m_i \times n_p}$ at a cost of $\mathcal{O}(m_i n_p r)$ operations and $\mathcal{O}(m_i n_p)$ bytes of memory for the result on the $i$th core.
To obtain the approximate solutions in the original coordinates, we apply all inverse data transformations, provided data transformations were used during preprocessing. 
These transformations are typically applied independently on each core, with a computational cost that varies depending on the specific transformation.
This cost, however, is usually low.
Finally, based on the postprocessing objective, the computed approximate solutions might be saved to disk or utilized for computing errors or other quantities.

\subsection{Algorithm summary} \label{subsec:dOpInf_summary}
In summary, dOpInf offers a complete distributed workflow for physics-based ROMs tailored to problems with extremely large state dimensions.
Key aspects of our approach include efficient data transformations and efficient dimensionality reduction without necessitating the computation of the POD basis. 
The primary computational steps involve standard linear algebra operations, such as matrix-matrix multiplications and eigenvalue solvers, which can be efficiently implemented using state-of-the-art scientific computing libraries.
This, and the fact that our formulation requires few communication steps, results in a scalable, distributed memory approach for learning physics-based ROMs.
These operations can be efficiently performed on CPUs as well as on specialized hardware such as graphical processing units or tensor processing units~\cite{Le22}. 

Elements of our distributed approach are transferable to other approaches such as DMD~\cite{Ku16, Sc10, Tu14} and quadratic manifolds~\cite{barnett2022quadratic, GWW23}, and also to parametric ROMs~\cite{BGW15}.
A summary for parametric ROM settings is as follows.
Let $\bm{\upmu} \in \mathbb{R}^{n_{u}}$ denote the $n_{u}$-dimensional vector of parameters of interest (e.g., comprising the mass flow, equivalence ratio and other relevant parameters in RDRE simulations).
Let $\bm{\upmu}_1, \bm{\upmu}_2, \ldots, \bm{\upmu}_z$ denote $z$ instances used for training and denote by $\mathbf{S}_{\alpha} \in \mathbb{R}^{n \times n_t}$ the corresponding snapshot matrices containing $n_t$ snapshots for $\alpha = 1, 2, \ldots, z$.
Following the parametric ROM survey in Ref.~\cite{BGW15}, two strategies for reducing the high-dimensional data are to use a global reduced basis using all $z$ snapshot matrices or separate reduced basis for each training parameter instance.
For either strategy Algorithm~\ref{alg:dOpInf_preprocessing} can be used for distributed data transformations and dimensionality reduction.
These approaches proceed by constructing reduced operators for each training parameter instance which can be done via Algorithm~\ref{alg:dOpInf_learning}.
Given a new parameter instance $\bm{\upmu}$, the corresponding reduced operators are typically obtained by interpolating the reduced operators for training; we refer to~\cite{BGW15} for more details on interpolation strategies.
Finally, the resulting ROM solution can be postprocessed in parallel analogously to Section~\ref{subsec:dOpInf_postprocessing}.

\section{Application to a large-scale real-world combustion scenario} \label{sec:results}
In this section, we assess the scalability and predictive capabilities of the proposed dOpInf algorithm in a complex three-dimensional unsteady RDRE scenario.
We focus on RDRE simulations as a representative example of a real-world application for which ROMs are needed to enable real-world engineering tasks that would be intractable with high-fidelity models.
Section~\ref{subsec:results_RDRE_overview} summarizes the considered RDRE scenario, followed by the description of the available high-fidelity dataset in Sec.~\ref{subsec:results_training}.
Section~\ref{subsec:results_scalability} presents our dOpInf scalability results using up to $2,048$ cores on the Frontera supercomputer.
Finally, we assess the prediction capabilities of dOpInf ROMs in the RDRE scenario under consideration in Sec.~\ref{subsec:results_predictions}.

\subsection{Overview of the considered rotating detonation rocket engine scenario} \label{subsec:results_RDRE_overview}
The RDRE concept injects fuel into an axially symmetric chamber, such as an annulus. 
When ignited under appropriate conditions, this process generates a system of spinning detonation waves~\cite{By06}. 
RDREs offer several advantages compared to conventional devices, including mechanical simplicity, which has led to active research in RDRE designs. 
Designing optimized RDREs, however, remains an open challenge due the difficulty in both simulation-based and experimental design space explorations. 
Large-eddy simulations (LES) of RDREs have advanced and can now provide valuable data that inform performance, stability, and realizability assessments of a given RDRE design~\cite{Ba20, Li19}. 
However, the large computational cost of LES, even on large supercomputers, makes it impractical for design optimization purposes---a single LES usually requires millions of core hours~\cite{Bat21}.
This highlights the need for scalable and predictive ROMs.

The physics of the considered RDRE scenario are modeled using the three-dimensional, reactive, viscous Navier-Stokes equations coupled with a skeletal chemistry mechanism (FFCMy-12) based on the FFCM model~\cite{FFCM1}.
This scenario is based on a design with $72$ discrete injector pairs.
The high-fidelity simulations for the full RDRE were performed using implicit LES via the AHFM (ALREST High-Fidelity Modeling) large-scale simulation code from the Air Force Research Laboratory (AFRL).
For more details about the code and typical setups, we refer the reader to \cite{Ba20, Li19}.
The high-fidelity LES employed a spatial discretization via a multi-block mesh comprising $136$ million spatial cells and a temporal discretization via an adaptive timestep $\Delta t \approx 10^{-9}$ seconds in the quasi-limit-cycle regime.
Each LES analysis requires approximately one million CPU-hours for one millisecond of simulated physical time on the DoD supercomputer Nautilus using $14,336$ cores across $112$ AMD EPYC Milan nodes.
This amounts to almost three days of CPU time, excluding queue waiting and postprocessing times.

\begin{figure}[htb!]
\centering
\includegraphics[width=1.0\textwidth]{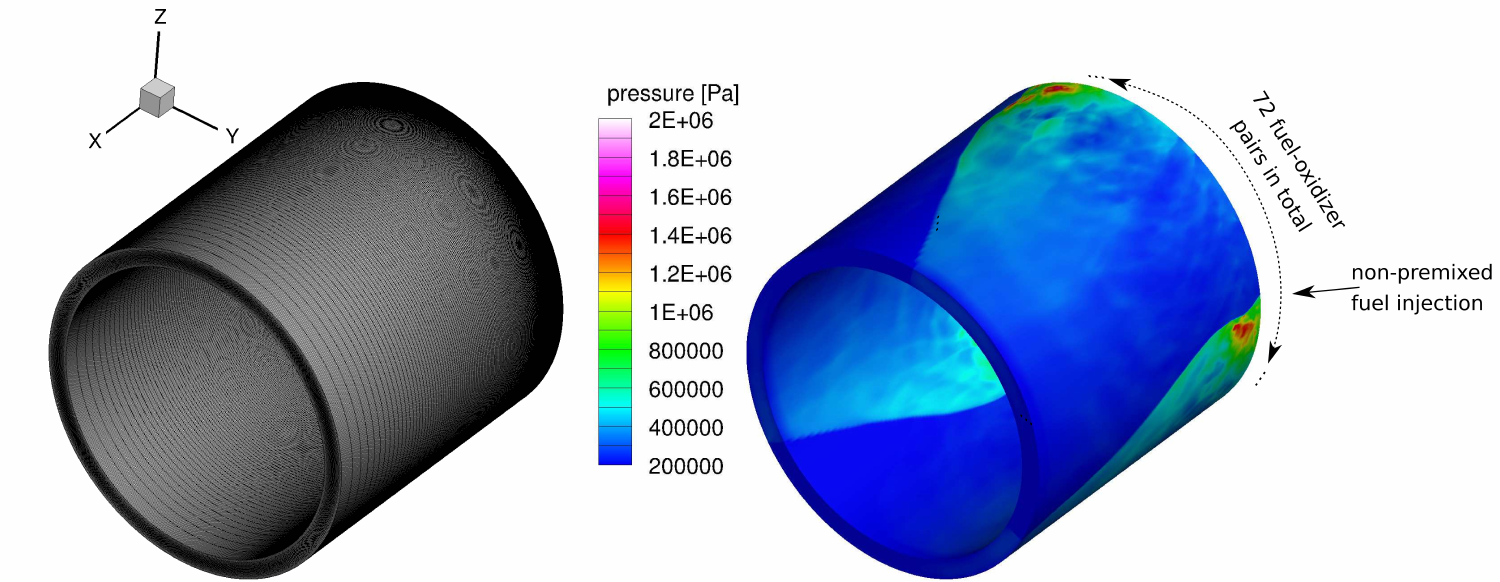}
\caption{Combustion chamber domain for the considered RDRE scenario. The left figure plots the structured mesh. The right figure plots an example pressure field showing the three dominant co-rotating waves.}
\label{fig:RDE_domain_and_pressure_ex}
\end{figure}

The input values of the mass flow-rate, $\dot{\textrm{m}} = 0.267 \ \textrm{kg} \cdot \textrm{s}^{-1}$, and equivalence ratio, $\Phi = 1.16$, lead to the formation of three dominant co-rotating waves and no secondary waves in the quasi-limit cycle regime.
Figure~\ref{fig:RDE_domain_and_pressure_ex} plots the combustion chamber of this RDRE, showing the mesh and an example pressure field with the three detonation waves.

We construct ROMs via dOpInf for the RDRE combustion chamber.
The three-dimensional computational domain spans from $0.05$ to $76.15$ millimeters in the $x$ direction and from $-37.50$ to $37.50$ millimeters in both $y$ and $z$ directions.
It has a fixed channel height of $4.44$ millimeters throughout.
The combustion chamber degrees of freedom were extracted from the full RDRE simulation results and interpolated onto a structured grid comprising $n_x = 4,204,200$ spatial degrees of freedom, with clustering of grid points at mid channel and closer toward the injector plane, as shown in Fig.~\ref{fig:RDE_domain_and_pressure_ex}.

\subsection{Acquiring the large-scale training dataset}~\label{subsec:results_training}
The governing equations are not in polynomial form.
We follow the work in~\cite{Sw20} and represent the transformed governing equations in terms of the following $m_s = 18$ physical state variables:    
\begin{equation} \label{eq:RDE_state_vars}
\begin{split}
    q = [ & 1/\rho \ \ p \ \  v_x \ \ v_y \ \ v_z \ \ w_{\mathrm{CH_4}} \ \ w_{\mathrm{O_2}} \ \ w_{\mathrm{CO_2}} \ \ w_{\mathrm{H_2O}} \ \ w_{\mathrm{CO}} \\ & w_{\mathrm{H_2}} \ \ w_{\mathrm{OH}} \ \ w_{\mathrm{CH_2O}} \ \ w_{\mathrm{CH_4}} \ \ w_{\mathrm{HO_2}} \ \ w_{\mathrm{H}} \ \ w_{\mathrm{O}} \ \ T]^\top,
\end{split}
\end{equation}
where $\rho \ [\mathrm{kg} \cdot \mathrm{m}^{-3}]$ is density (and $1/\rho \ [\mathrm{m}^{3} \cdot \mathrm{kg}^{-1}]$ is specific volume), $p \ [\mathrm{Pa}]$ is pressure, $v_x, v_y, v_z \ [\mathrm{m} \cdot \mathrm{s}^{-1}]$ are the velocity components, $w_{\cdot}$ are the mass fractions of all chemical species ($12$ in total here), and $T \ [\mathrm{K}]$ is temperature.
The $m_s = 18$ variables in Eq.~\eqref{eq:RDE_state_vars} make most terms in the lifted governing equations linear or quadratic.
The dimensionality of the transformed snapshots is therefore $m = 18 \times 4,204,200 = 75,675,600$.

The computational cost and size of the resulting simulation datasets limit the number of time instants for which the high-fidelity LES solutions can be saved to disk.
We have a total of $3,805$ downsampled snapshots from the high-fidelity LES over the time horizon $[5.0000, 5.3804]$ milliseconds (the down-sampling factor was about $100$), which represents one of the largest RDRE datasets used for ROM development.
The $3,805$ downsampled snapshots correspond to about three full quasi-cycles of the three-wave system.
Because of the severe downsampling of the snapshots, we employ the discrete formulation of OpInf~\cite{Fa23} (cf.~Section~\ref{subsec:dOpInf_learning}).
Wave cycles in the quasi-steady state of RDRE simulations usually exhibit near-periodicity but also cycle-to-cycle variations. 
This lack of self-similarity can present a considerable challenge for accurate modeling.
Simulating the full $5.3804$ milliseconds of physical time required about $5.4$ million core hours on $14,336$ compute cores, which means that the $0.3804$ milliseconds spanning the target time interval necessitated about $380,000$ core hours, corresponding to a simulation time of $27$ hours. 
We use the first $n_t = 2,536$ snapshots (i.e., two thirds of the total number of available snapshots) for training dOpInf ROMs, corresponding to two quasi-cycles of the three-wave system.
The remaining $1,269$ snapshots are used to assess the prediction capabilities of the dOpInf ROM beyond training.
We note that constructing ROMs for RDRE predictions beyond a training time horizon is an essential step for assessing the prediction capabilities of ROMs in these complex applications as well as for enabling predictions over time horizons that would be computationally too expensive to simulate using high-fidelity LES---recall that one millisecond of simulated physical time via LES can require one million core hours on a supercomputer.
An initial study showcasing the potential of OpInf (using the standard, serial formulation) for constructing parametric data-driven ROMs for three-dimensional unsteady RDREs can be found in Ref.~\cite{Fa23}.
As noted in Section~\ref{subsec:dOpInf_summary}, elements of our dOpInf workflow can be also transferred to the parametric setting.

The large dimension $m = 75,675,600$ of the $n_t = 2,536$ training snapshots means that their total size amounts to $1.4$ TB in double precision. 
As such, performing data transformation and dimensionality reduction to construct ROMs would be computationally expensive and limited by memory storage and bandwidth.
Our distributed formulation, in contrast, offers a scalable workflow for constructing physics-based ROMs on distributed memory machines.

\subsection{Scalability results on the Frontera supercomputer}~\label{subsec:results_scalability}
We first demonstrate the scalability of the proposed dOpInf approach, from loading the training data in parallel to computing the reduced solution over the target time horizon $[5.0000, 5.3804]$ (i.e., using Algorithms~\ref{alg:dOpInf_preprocessing} and~\ref{alg:dOpInf_learning}). 
The number of training  snapshots, $n_t = 2,536$, permits a maximum reduced dimension $r=68$ for a quadratic dOpInf ROM since it sets the maximum number of operator coefficients that can be learned via the regularized regression problem~\eqref{eq:OpInf_reg}.
In all our subsequent experiments we construct dOpInf ROMs with reduced dimension $r=68$.
The number of training snapshots (i.e., $n_t = 2,536$) and the target time horizon (i.e., $[5.0000, 5.3804]$ milliseconds) are the same in all experiments.
Moreover, the trial time horizon for the optimal regularization parameter search in Algorithm~\ref{alg:dOpInf_learning} is the same as the target horizon.

Our distributed memory implementation was done using \texttt{Python} (\texttt{Python 3.7} on Frontera) and the Message Passing Interface (MPI) bindings provided by the \texttt{mpi4py} library.
We leveraged TACC's scientific software stack to maximize performance\footnote{\url{https://docs.tacc.utexas.edu/hpc/frontera/\#ml}}.
\texttt{mpi4py} uses the Intel MPI implementation tailored to the underlying hardware, thus ensuring fast MPI operations.
The dense linear algebra operations are performed using \texttt{numpy}'s linear algebra functions \texttt{matmul} (for matrix-matrix multiplications) and \texttt{scipy}'s \texttt{eigsh} function for partial eigendecompositions of real symmetric matrices.
\texttt{numpy} and \texttt{scipy} are linked to the Intel Math Kernel Library (MKL) and compiled to take advantage of CORE-AVX512 instruction set architecture relevant to Frontera nodes. 
MKL provides access to optimized versions of the low-level BLAS and LAPACK libraries.

To obtain a comprehensive assessment of our algorithm's scalability, we use up to $2,048$ cores distributed over Cascade Lake (CLX) compute nodes on the Frontera supercomputer. 
Each CLX node provides $56$ cores across two sockets and $192$ GB of main memory. 
In assessing our algorithm's scalability, we average the CPU times over five measurements for each core count to mitigate the effects of fluctuations.

The RDRE training dataset was saved in HDF5 format and transferred on the \texttt{scratch1} partition on the Frontera supercomputer.
This employs a Lustre file system\footnote{\url{https://docs.tacc.utexas.edu/hpc/frontera/\#files}}, typically used on large-scale computing systems.
To mitigate file system bottlenecks when using a large number of compute cores, we distributed the $1.4$ TB training dataset across multiple HDF5 files, which represents a common approach for handling large datasets.
Each file contains data for the $m_s = 18$ transformed physical state variables within a specific subdomain of the full computational domain.
For a comprehensive overview of Frontera's file system and recommended best practices, we refer the reader to Ref.~\cite{frontera_IO}.
For loading the training dataset in parallel in our implementation, we used the \texttt{h5py} library which is linked to the parallel HDF5 (phdf5) library on Frontera.
We note that data reading is independent of the other dOpInf algorithmic steps. 
Users are thus responsible for ensuring that their training data are stored in a format conducive to efficient I/O operations.

We begin with the results for strong scaling.
The goal is to assess the speed-ups obtained with the dOpInf algorithm, from loading the data in parallel to computing the reduced solution using the dOpInf ROM over the target time horizon $[5.0000, 5.3804]$ milliseconds.
We keep the problem size fixed, that is, we employ the entire snapshot data of size $75,675,600 \times 2,536$ to construct a predictive dOpInf ROM with reduced dimension $r=68$ using our distributed workflow, and increase the number of compute cores.
Here, we use $p \in \{ 32; 64; 128; 256; 512; 1,024; 2,048\}$.
The regularization hyperparameter search is conducted using a grid of size $32 \times 64$  (amounting to $B = 2,048$ candidate pairs in total).
The speed-up is computed as $T(32)/T(p)$.
Here, $T(p)$ denotes the average CPU time using $p$ cores measured on the compute core that contains the optimal regularization hyperparameters.
The ideal speed-up is $p/32$.

\begin{figure}[htb!]
\centering
\includegraphics[width=\textwidth]{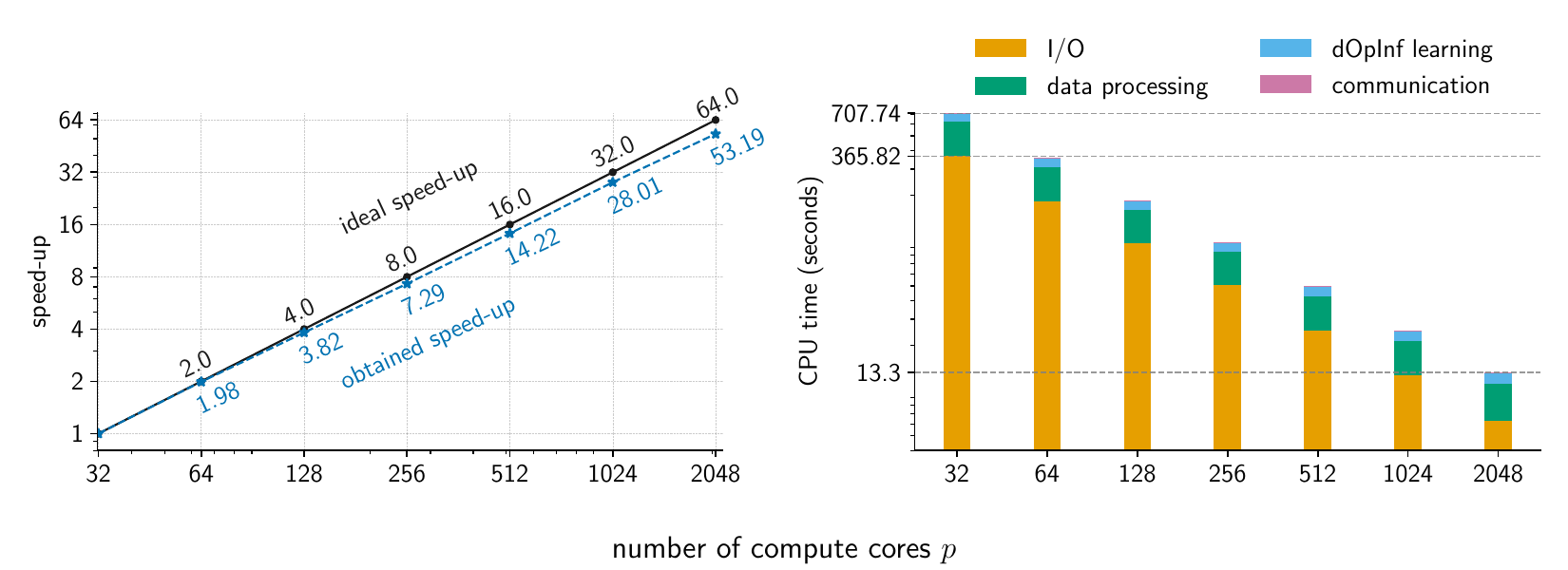}
\caption{Strong scaling results using between $p=32$ and $p=2,048$ cores on Frontera. The left plot shows the speed-ups and the right plot shows the percentages of the total CPU time corresponding to data loading, all data processing computations, learning the reduced operators with dOpInf, and communication overhead.}
\label{fig:RDE_strong_scaling}
\end{figure}
\begin{figure}[htb!]
\centering
\includegraphics[width=\textwidth]{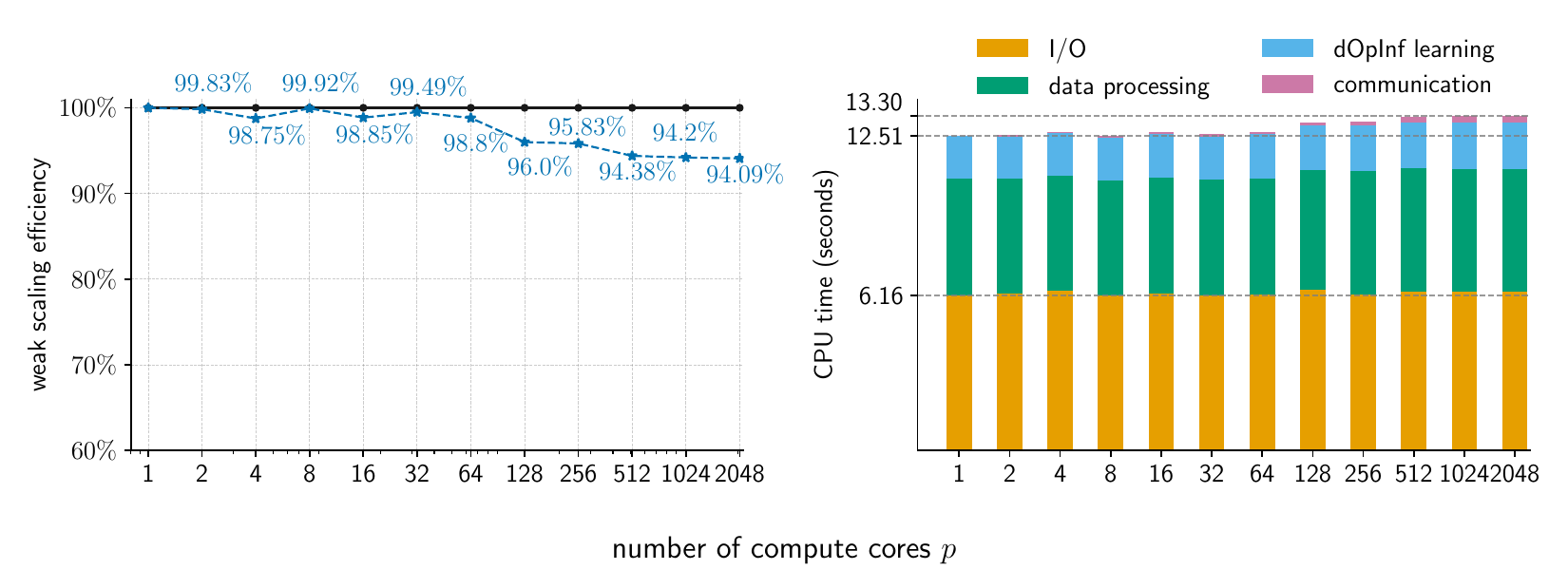}
\caption{Weak scaling results using between $p=1$ and $p=2,048$ cores on Frontera. The left plot shows the obtained efficiency and the right plot shows the corresponding CPU times for data loading, all data processing computations, learning the reduced operators with dOpInf, and communication overhead.}
\label{fig:RDE_weak_scaling}
\end{figure}

Figure~\ref{fig:RDE_strong_scaling} plots the results.
The left plot shows the speed-ups relative to $p=32$ cores and the right plot shows the corresponding CPU time distribution into loading the snapshot data (I/O), all data processing computations in Algorithm~\ref{alg:dOpInf_preprocessing}, learning the reduced operators via Algorithm~\ref{alg:dOpInf_learning}, and the communication overhead in both Algorithms~\ref{alg:dOpInf_preprocessing} and~\ref{alg:dOpInf_learning}.
The total CPU time decreases from $707.74 \pm 8.94$ seconds for $p=32$ cores down to $13.30 \pm 0.04$ seconds for $p=2,048$ compute cores. 
Our results closely match the ideal speed-ups, indicating excellent strong scalability.
I/O and data processing computations drive the CPU times, but their excellent scalability coupled with a minimal communication overhead leads to the obtained near-ideal speed-ups.
Because of this, constructing the dOpInf ROM with dimension $r=68$ took only $13.30 \pm 0.04$ seconds using $p = 2,048$ cores. 
This breaks down to $6.30 \pm 0.04$ seconds for I/O, $4.89 \pm 0.01$ seconds for all data processing computations, $1.84 \pm 0.01$ seconds for dOpInf learning, and $0.26 \pm 0.01$ seconds for communication.
Therefore, in contrast to many existing data-driven reduced and surrogate modeling approaches that require substantial CPU times, our method enables a rapid and efficient construction of physics-based ROMs in problems with large datasets and a large snapshot dimension.

We next assess the weak scaling efficiency of Algorithms~\ref{alg:dOpInf_preprocessing} and~\ref{alg:dOpInf_learning}.
In contrast to the speed-up analysis where we fixed the problem size and scaled the number of cores, here we fix the problem size per core. 
Our goal is to evaluate how effectively our distributed workflow utilizes resources as both the problem size (i.e., the size of the snapshot matrix and the resolution of the grid used for searching the optimal regularization parameter pair) and the number of cores increase, which provides more insights into its scalability for datasets with extremely large state dimension.
We vary the number of cores from $p=1$ and $p=2,048$, increasing in powers of two.
The largest core count, $p=2,048$, is used to construct the dOpInf ROM using the full training dataset and a grid of size $32 \times 64$ for the regularization hyperparameter search.
Note that this setup is the same as the one used for the strong scaling analysis for $p=2,048$ cores.
To maintain a consistent problem size per core for configurations with fewer than $2,048$ cores, we adjust the problem size by downsampling the full snapshot dimension and reducing the grid size for the hyperparameter search such that we have $B = p$ candidate pairs for all considered $p$ core counts.  
The computed weak scaling efficiency is $T(1)/T(p)$, where $T(1)$ corresponds to the serial implementation of Algorithms~\ref{alg:dOpInf_preprocessing} and~\ref{alg:dOpInf_learning} using a downsampled snapshot matrix of size $36,951\times 2,536$ (i.e., with a downsampled snapshot dimension $75,675,600/2,048$), averaged over five runs.
We obtain $T(1) = 12.51 \pm 0.20$ seconds ($6.16 \pm 0.15$ seconds for I/O, $4.66\pm 0.06$ for all data processing computations, and $1.68 \pm 0.02$ seconds for dOpInf learning).

A good weak scaling performance would be obtained if our algorithm maintained a similar CPU time for I/O and all computations, and a low communication overhead as $p$ increases. 
Figure~\ref{fig:RDE_weak_scaling} plots the results. 
Our algorithm indeed demonstrates an excellent weak scaling efficiency that exceeds $94\%$ for all values of $p$, as shown on the left.
The right plot in Fig.~\ref{fig:RDE_weak_scaling} shows that our algorithm maintains a similar CPU time for data loading and all computations, whereas the communication overhead increases only slowly with $p$.

Two important elements for the obtained scalability results were the efficient parallel loading of the large training dataset and the low computational cost of the sequential steps in Algorithm~\ref{alg:dOpInf_preprocessing} (cf.~Remarks~\ref{remark:collective_reduction_Algo_1} and~\ref{remark:scalability_Algo_1}).
Furthermore, the dOpInf efficiency, evident in its ability to complete the entire workflow (parallel data loading, data transformations, dimensionality reduction, and ROM construction via a grid search over $B = 2,048$ candidate pairs) in an average time of just $13.30$ seconds on $p=2,048$ cores was facilitated by the powerful scientific software stack on Frontera.

\subsection{Predictions beyond the training time horizon} ~\label{subsec:results_predictions}
We now assess the potential of the proposed distributed memory model reduction workflow for constructing predictive structure-preserving, physics-based ROMs for the real-world RDRE scenario under consideration.
We use the ROM results obtained for the scalability studies above for $p=2,048$ cores.
Recall that the average CPU time for constructing the dOpInf ROM with reduced dimension $r=68$ amounted to $13.30$ seconds in total.

\begin{remark} \label{remark:ROM_adequency}
Due to the complexity of RDRE simulations, we prioritize ROMs that accurately capture large-scale features, including detonation and shockwave fronts. 
These ROMs should preserve essential time-averaged chamber quantities, wave propagation characteristics, and key engineering quantities like pressure, axial velocity/temperature, and fuel/oxidizer mass fractions. 
Achieving a ROM that accurately represents these large-scale features provides a practically useful capability for design exploration and analysis tasks, even if it does not capture all small-scale features.
\end{remark}

\begin{figure}[htb!]
\centering
\includegraphics[width=0.8\textwidth]{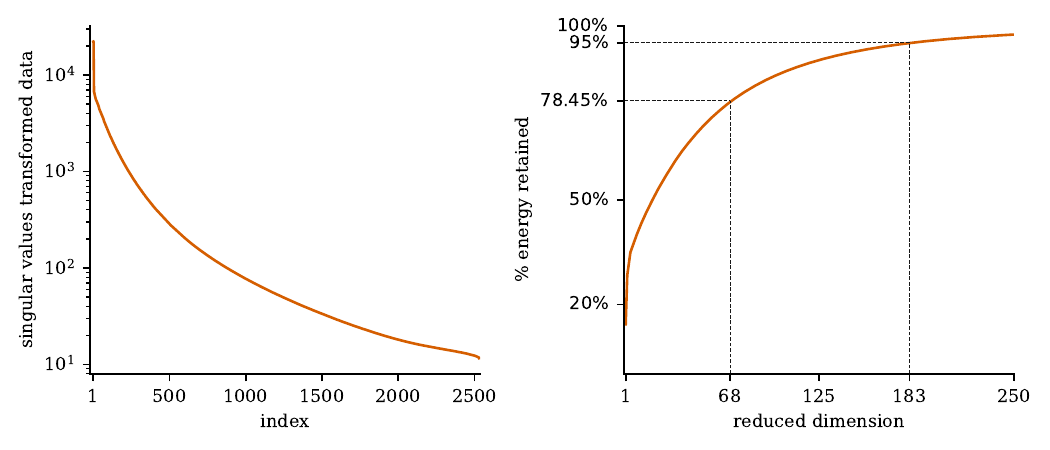}
\caption{The left figure plots the POD singular values of the transformed snapshots. The right figure plots the corresponding retained energy.}\label{fig:RDE_svals_and_ret_energy}
\end{figure}

We start by analyzing the decay of the POD singular values.
Figure~\ref{fig:RDE_svals_and_ret_energy} plots the POD singular values on the left and corresponding retained energy on the right.
A slow decay of the singular values is expected due to the complex dynamics in this scenario. 
In addition, the corresponding POD basis is global, representing all $m_s = 18$ transformed variables in~\eqref{eq:RDE_state_vars}, characterized by heterogeneous dynamics.
Retaining $95\%$ of the total energy requires a large reduced dimension, $r = 183$. 
In contrast, $r = 68$ retains $78.45\%$ of the total energy.
We will nonetheless show that the dOpInf ROM with reduced dimension $r=68$ produces predictions that meet the requirements specified in Remark~\ref{remark:ROM_adequency}. 

We then assess how well the ROM solution meets the accuracy criteria presented in Remark~\ref{remark:ROM_adequency}.
We utilize the ROM solution to extract one-dimensional radial profiles for pressure, temperature, and fuel ($\mathrm{CH_4}$) mass fraction at three representative locations near the mid-channel of the combustion chamber.
Axially, the first location is close to the injectors, the second location is further way from the injectors but still within the detonation region, and the third location is downstream of the detonation zone.
For a more comprehensive evaluation, we also compute the full pressure field. 
This enables us to assess the ROM's ability to represent the large-scale features of the three detonation waves in the full computational domain.
To map the ROM solution back to the original coordinates, we compute the components of the global POD basis of rank $r=68$ via Eq.~\eqref{eq:MoS_POD_basis_comp}, which requires an average of $0.24$ seconds on $p=2,048$ cores.

\begin{figure}[htb!]
\centering
\includegraphics[width=\textwidth]{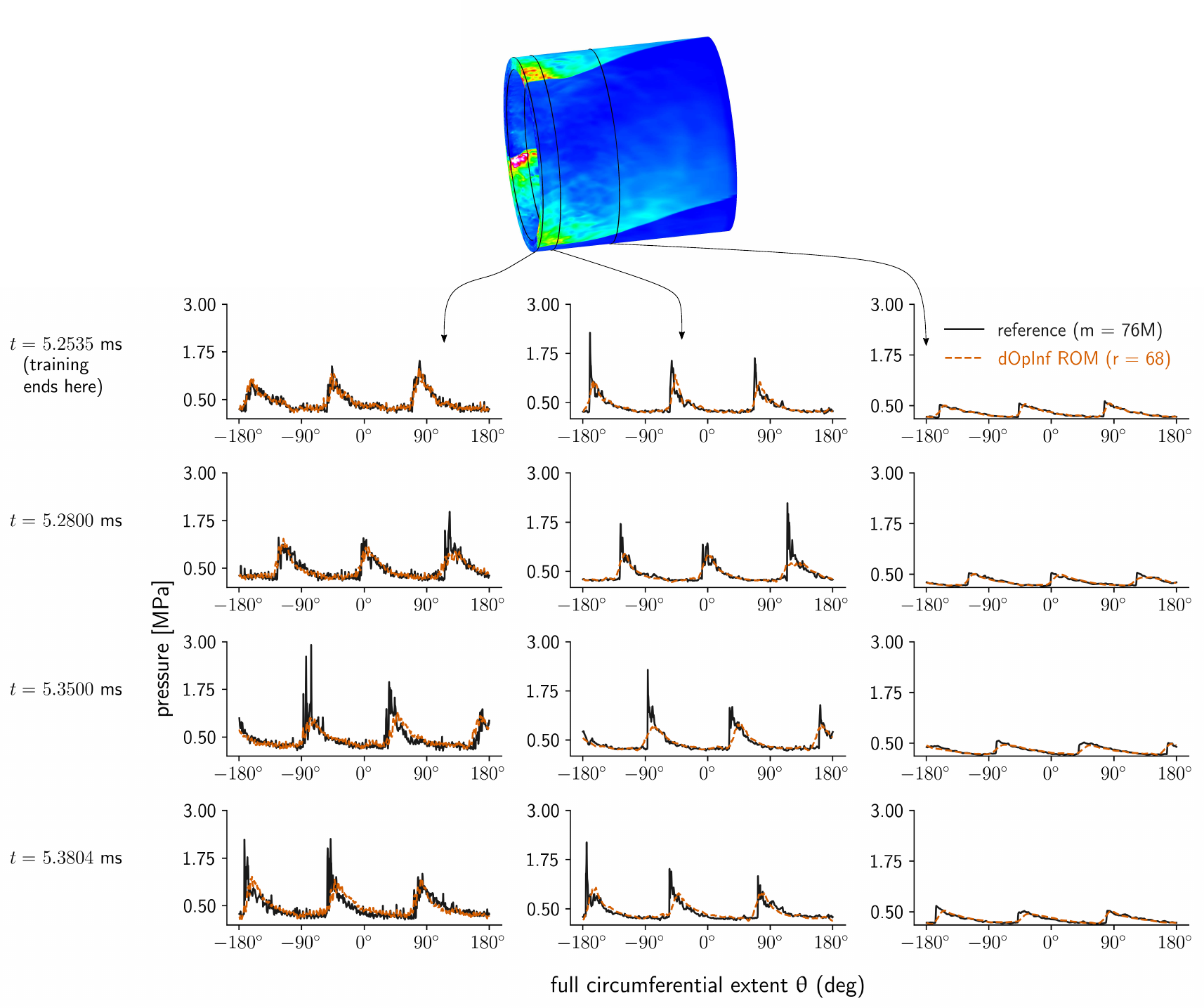}
\caption{One-dimensional circumferential profiles for pressure. The columns plot the results at three representative locations close to the mid-channel. The rows plot the profiles at four representative time instants.}
\label{fig:RDE_1D_pressure_profiles}
\end{figure}

\begin{figure}[htb!]
\centering
\includegraphics[width=\textwidth]{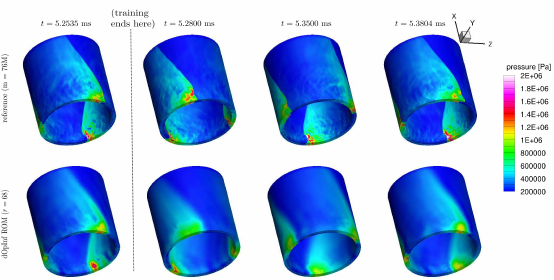}
\caption{Pressure fields at four time instants. The top row plots the reference solutions. The bottom row plots the approximate solution obtained with our reduced model.}
\label{fig:RDE_pressure_isosurfaces}
\end{figure}

\begin{figure}[htb!]
\centering
\includegraphics[width=\textwidth]{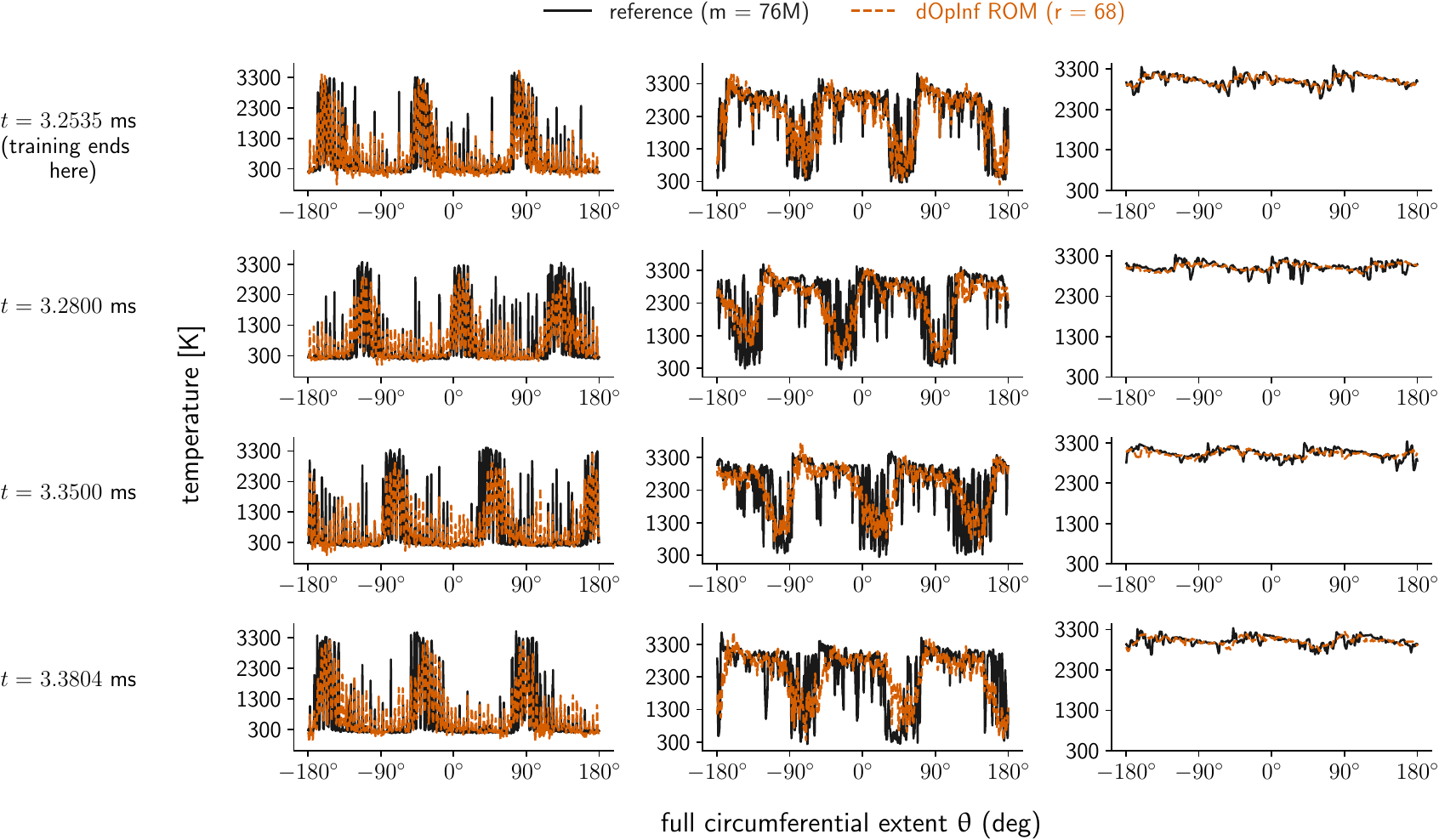}
\caption{One-dimensional circumferential profiles for temperature. The columns plot the results at three representative locations close to the mid-channel. The rows plot the profiles at four representative time instants.}
\label{fig:RDE_1D_temperature_profiles}
\end{figure}

\begin{figure}[htb!]
\centering
\includegraphics[width=\textwidth]{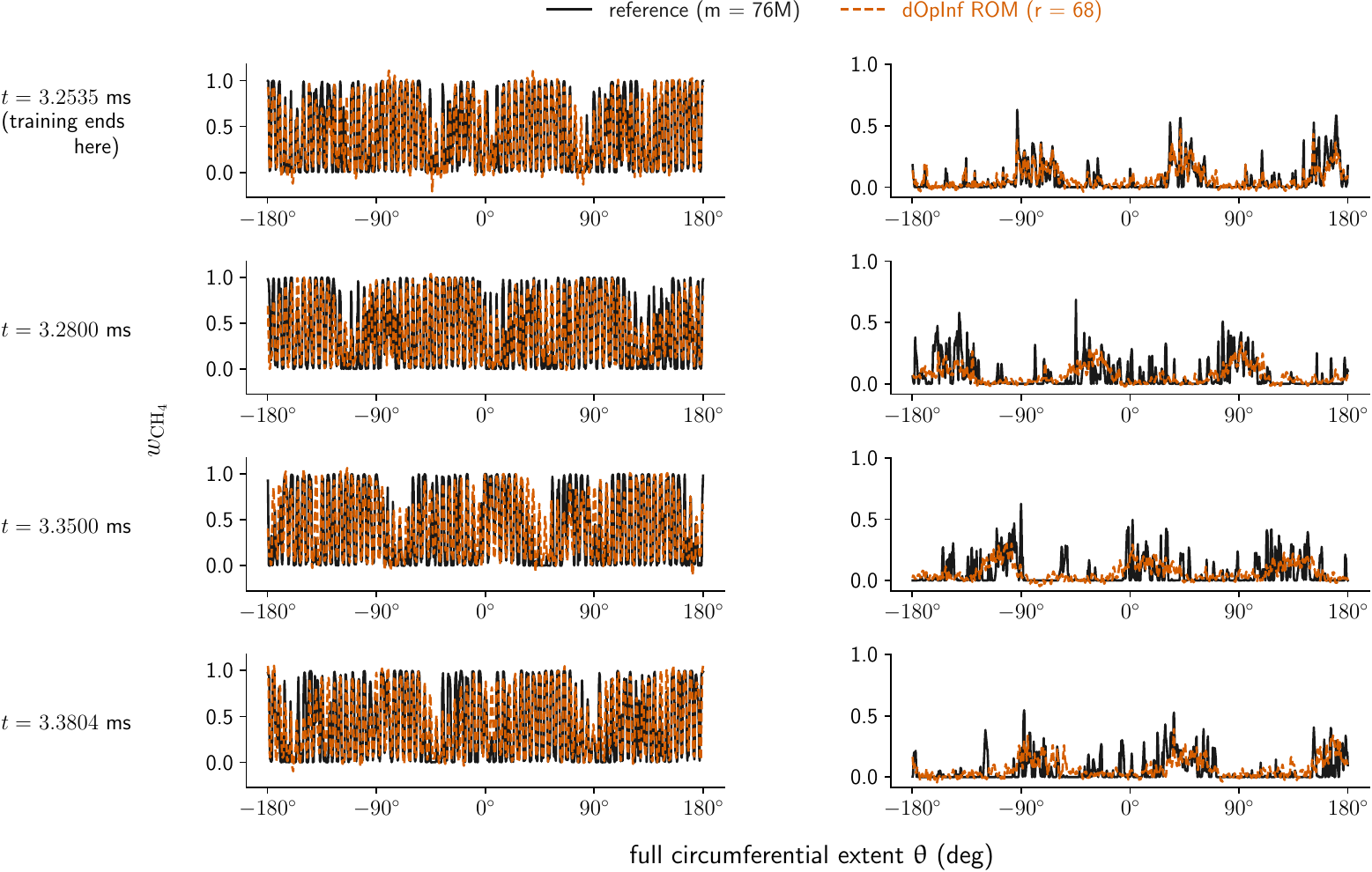}
\caption{One-dimensional circumferential profiles for fuel ($\mathrm{CH_4}$) mass fraction. The columns plot the results at two representative locations close to the mid-channel. The rows plot the profiles at four representative time instants.}
\label{fig:RDE_1D_CH4_profiles}
\end{figure}

We present first the approximate solutions for pressure.
Figure~\ref{fig:RDE_1D_pressure_profiles} plots the radial profiles at four representative time instants: the last time instant in the training horizon, $t = 5.2535$ ms, and three time instances in the prediction horizon, $t = 5.2800$ ms (the $265$th), $t = 5.3500$ ms (the $965$th), and $t = 5.3804$ ms (the last, i.e., $1269$th time instant in the prediction horizon).
The obtained profiles align well with the corresponding reference results: the ROM approximate solutions capture larger-scale features, including the number of waves, wave fronts, and also the lower amplitudes of the reference solutions, but, as expected, they do not capture the abruptly changing pressure spikes in the reference solution.
Nevertheless, pressure spikes are not considered a critical metric in this context, as they can exhibit significant variability among high-fidelity simulations.
Figure~\ref{fig:RDE_pressure_isosurfaces} plots the corresponding full fields, which shows that the ROM solution adequately captures the key characteristics of the shock waves in the full computational domain.

Figures~\ref{fig:RDE_1D_temperature_profiles} and~\ref{fig:RDE_1D_CH4_profiles} plot the radial profiles for temperature and fuel mass fraction, respectively, at the same time instants as for the pressure profiles in Figure~\ref{fig:RDE_1D_pressure_profiles}.
Note that Figure~\ref{fig:RDE_1D_CH4_profiles} plots the profiles at only the first two locations since the fuel mass fraction values at the third location are close to zero.
The non-smoothness of these profiles as well as the variations between the three quasi-cycles further evidentiate the complexity of the dynamics in the RDRE scenario under consideration.
However, our dOpInf ROM captures well the larger-scale features.

Finally, the average evaluation time of the resulting dOpInf ROM is $1.09$ seconds on one CLX core.
This represents a computational cost reduction of $90,000$ compared to the CPU time of the high-
fidelity simulation over the target time horizon.

\section{Conclusion} \label{sec:conclusions}
Recent advancements in HPC simulations of complex real-world problems necessitate the development of innovative, parallelizable data-driven model reduction techniques tailored to modern HPC architectures. 
This paper demonstrated the power of integrating HPC into data-driven reduced modeling. 
The proposed distributed Operator Inference algorithm allows a fast and scalable processing of extremely large datasets, and the construction of predictive physics-based reduced models that approximate the dynamics underlying these datasets. 
These capabilities unlock new possibilities for computationally expensive tasks like design optimization, which would otherwise be intractable using high-fidelity models.
These developments hold promise for a wide range of fields, including rocket propulsion and the assessment of turbulent transport in fusion devices.
An implementation of the distributed Operator Inference algorithm, including a detailed tutorial, is available at \url{https://github.com/ionutfarcas/distributed_Operator_Inference}.

The proposed distributed algorithm inherits the limitations of standard Operator Inference, namely, the difficulty of effectively constructing reduced models for problems characterized by slowly decaying Kolmogorov $n$-widths.
This can be addressed, for example, via quadratic manifolds~\cite{barnett2022quadratic, GWW23}, which can be  extended to a distributed formulation following a similar approach to that presented here.

Domain decomposition also holds promise for parallel processing of large datasets~\cite{Fa24}.
However, methods like Operator Inference struggle with a high number of subdomains due to hyperparameter tuning complexity. 
Finding a balance between subdomain usage and accuracy, potentially through integration with distributed methods, could unlock significant efficiency gains.

\section*{Acknowledgements}
This work was supported in part by AFRL Grant FA9300-22-1-0001 and the Air Force Center of Excellence on Multifidelity Modeling of Rocket Combustor Dynamics under grant FA9550-17-1-0195.
The authors gratefully acknowledge Jonathan Hoy who helped with transferring the high-fidelity dataset to TACC.
The authors also gratefully acknowledge the compute and data resources provided by the  Texas Advanced Computing Center at The University of Texas at Austin \url{https://www.tacc.utexas.edu} and the DoD High Performance Computing Modernization Program (HPCMP).
The views expressed are those of the author and do not necessarily reflect the official policy or position of the Department of the Air Force, the Department of Defense, or the U.S. government.

Distribution Statement A: Approved for Public Release; Distribution is Unlimited. PA$\#$ AFRL-2024-111

\bibliographystyle{elsarticle-harv}
\bibliography{distributed_OpInf.bib}

\end{document}